\documentclass{siamltex}%
\usepackage{amssymb}
\usepackage{amsmath}
\usepackage{amsfonts}
\usepackage{sw20siam}
\usepackage{graphicx}%
\providecommand{\U}[1]{\protect\rule{.1in}{.1in}}

\numberwithin{equation}{section}

\begin{document}

\title{Energy Dissipation in the Smagorinsky Model of Turbulence}
\author{William Layton\thanks{\bigskip Department of Mathematics, University of
Pittsburgh, Pittsburgh, PA 15260, USA, wjl@pitt.edu,
www.math.pitt.edu/\symbol{126}wjl . Partially supported by NSF grant DMS
1216465 and 1522267.}}
\date{May 2016}
\maketitle

\begin{abstract}
The Smagorinsky model, unmodified, is often reported to severely overdiffuse
flows. Previous\ estimates of\ the energy dissipation rate of the Smagorinsky
model for \textit{shear flows} are that%
\[
\left\langle \varepsilon_{S}\right\rangle \simeq\lbrack1+C_{S}^{2}\left(
\frac{\delta}{L}\right)  ^{2}(1+Re^{2})]\frac{U^{3}}{L}%
\]
\ reflecting a blow\ up\ of model energy dissipation as $Re\rightarrow\infty$.
This blow up is consistent with the numerical evidence and leads to the
question: \textit{Is the over dissipation due to the influence of the
turbulent viscosity in boundary lauyers alone or is its action on small scales
generated by the nonlinearity through the cascade also a contributor?} This
report develops model dissipation estimates for body force driven flow under
periodic boundary conditions (and thus only with nonlinearity generated small
scales). It is proven that the model's time averaged energy dissipation rate,
$<\varepsilon_{S}>$, satisfies%
\[
<\varepsilon_{S}>\leq3\frac{U^{3}}{L}+\frac{3}{8}\operatorname{Re}^{-1}%
\frac{U^{3}}{L}+C_{S}\left(  \frac{\delta}{L}\right)  ^{2}\frac{U^{3}}%
{L}\text{,\ }%
\]
where $U,L$ are global velocity and length scales and $C_{S}\simeq0.1$,
$\delta<<1$\ are the standard model\ parameters. Since this estimate is
consistent with that observed for the NSE, it establishes that,
\textit{without boundary layers,\ the Smagorinsky\ model does\ not\ over
dissipate.}

\end{abstract}

\textbf{Key words}. energy dissipation rate, Smagorinsky model, large eddy
simulation, turbulence


\section{Introduction}

\begin{center}
\textit{This is an expanded version of a report with a similar title.}
\end{center}

The Smagorinsky model\footnote{In its most precise realization, the term
$\nabla\cdot\left(  \left(  C_{S}\delta\right)  ^{2}|\nabla u|\nabla u\right)
$ is replaced by $\nabla\cdot\left(  2\left(  C_{S}\delta\right)  ^{2}%
|\nabla^{s}u|\nabla^{s}u\right)  $\ where $\nabla^{s}u:=(\nabla u+\nabla
u^{T})/2$ is the deformation tensor. The analysis herein holds by the same
argument.}, from \cite{S63} and used for of turbulent flow, e.g.,
\cite{BIL06}, \cite{Sagaut}, \cite{VolkersBook}, \cite{Mus96}, \cite{P94},
\cite{Geu97}, is%
\begin{gather}
u_{t}+u\cdot\nabla u-\nu\triangle u-\nabla\cdot\left(  \left(  C_{S}%
\delta\right)  ^{2}|\nabla u|\nabla u\right)  +\nabla
p=f\,(x)\,\,\label{eq:SmagorinskyModel}\\
\text{and }\nabla\cdot u=0\,\,.\nonumber
\end{gather}
It is mathematically equivalent to\ the Ladyzhenskaya model\footnote{In the
Smagorinsky model one has the viscous terms $-\nu\triangle u-\nabla
\cdot\left(  \left(  C_{S}\delta\right)  ^{2}|\nabla u|\nabla u\right)  $
while in the Ladyzhenskaya model the corresponding terms are $-\nabla
\cdot\left(  \sqrt{\nu^{2}+\left(  C_{S}\delta\right)  ^{4}|\nabla u|^{2}%
}\nabla u\right)  .$ For both, the most precise realization replaces the
velocity gradient with the deformation tensor.} \cite{S63}, \cite{Lad2},
\cite{DG91}, \cite{P92} and the von Neumann Richtmyer artificial viscosity for
shocks \cite{vNR}. \ In (\ref{eq:SmagorinskyModel}), $\nu$ is the kinematic
viscosity, $\delta<<1$ is a model length scale, the Reynolds number is
\[
Re=\frac{LU}{\nu}%
\]
where $U,L$ denote global velocity and length scales given by
(\ref{eq:ULscales}) below and $C_{S}$ is a model parameter. See \cite{S84} for
ranges of values around the value $C_{S}\simeq0.1$ determined by Lilly
\cite{L67}. Experience with the model, e.g., \cite{Sagaut}, \cite{BIL06},
strongly suggests it over dissipates, often severely (see Section 3 for some
fixes). \ Estimates of model energy dissipation rates for shear flows in
\cite{L02} are consistent with this computational experience. Perhaps
surprisingly, herein we show that the time averaged\ energy dissipation rate
for (\ref{eq:SmagorinskyModel})\ balances the energy input rate\footnote{The
energy input rate at the large scales is $U^{3}/L$. Briefly, the kinetic
energy of the large scales scales with dimensions $U^{2}$. The "\textit{rate}"
has dimensions \textit{1/time}. A large scale quantity with this dimensions is
formed by $U/L$ which is the turn over time for the large eddies, i.e., the
time iit takes a large eddy with velocity $U$ to travel a distance $L$. Thus
the "\textit{rate of energy input}" has dimensions $U^{3}/L$.}, $U^{3}/L$.
Thus \textit{the Smagorinsky model does not\ over dissipate energy for body
force driven turbulence in a periodic box}. In other words, the observed over
dissipation of the model is not due to the model's action on small scales
generated by the nonlinear term in the turbulent cascade but rather it is
\textit{due to the action of the model viscosity in boundary layers.}

\textit{How does the fluid velocity produced from a numerical simulation on a
fixed grid communicate with molecular viscosity? }- J. Smagorinsky 1960

The motivation for the model can be described in a very general way as
follows. For high Reynolds numbers, dissipation occurs non-negligibly only at
very small scales, far smaller than typical meshes. The balance between energy
input at the largest scales and energy dissipation at the smallest is a
critical selection mechanism for determining statistics of turbulent flows.
Joseph Smagorinsky was concerned with geophysical flow simulations and soon
asked the above question. Somehow, once a mesh is selected, to get accurate
simulations extra dissipative terms must be introduced to model the effect of
the unresolved fluctuations (smaller than the mesh width) upon the resolved
velocity (representable on the mesh). Thus the extra term was in
(\ref{eq:SmagorinskyModel}) was introduced where the length scale $\delta$ was
intended to reflect the underlying physical mesh length.

Let $\Omega=(0,L_{\Omega})^{3}$ denote the periodic box in $3d$ and impose
periodic (with zero mean) conditions
\begin{align}
u(x+L_{\Omega}e_{j},t)  &  =u(x,t)\quad j=1,2,3\text{ and }\label{eqBC}\\
\int_{\Omega}\phi dx  &  =0\,\text{\ for }\,\phi=u,\,u_{0},\,f,\,p.\nonumber
\end{align}
The data $u_{0}(x),f(x)$\ are smooth, $L_{\Omega}$-periodic, have\ zero mean
and satisfy%
\[
\nabla\cdot u_{0}=0\text{ , and }\nabla\cdot f=0.
\]
The model energy dissipation rate from (\ref{eq:EnergyEquality}) below is%
\[
\varepsilon_{S}(u):=\int_{\Omega}\frac{\nu}{|\Omega|}|\nabla u(x,t)|^{2}%
+\frac{(C_{S}\delta)^{2}}{|\Omega|}|\nabla u(x,t)|^{3}dx.
\]
The\ long time average of a function $\phi(t)$\ is defined, following
\cite{DF02}, \cite{DG95} by%
\[
\left\langle \phi\right\rangle :=\lim\sup_{T\rightarrow\infty}\frac{1}{T}%
\int_{0}^{T}\phi(t)dt.
\]
We show herein that $\left\langle \varepsilon_{S}\right\rangle $\ balances the
energy input rate, $U^{3}/L$. This estimate is consistent as $Re\rightarrow
\infty,\delta\rightarrow0$ with both phenomenology, e.g., \cite{Frisch},
\cite{Pope}, \cite{Mus96}, and the rate proven for the Navier-Stokes equations
in \cite{CD92}, \cite{CKG01}, \cite{Wang97}, \cite{F97}, \cite{DF02} and
\cite{CDP06}. The weak $Re$ dependence in the second term (that vanishes as
$Re$ $\rightarrow\infty$) is consistent with the recent results in
\cite{MBYL15} derived through structure function theories of turbulence.

\begin{theorem}
Suppose the data $f(x)$ and $u_{0}(x)$ are smooth, divergence free, periodic
with zero mean functions. Then%
\[
\left\langle \varepsilon_{S}(u)\right\rangle \leq3\frac{U^{3}}{L}+\frac{3}%
{8}\operatorname{Re}^{-1}\frac{U^{3}}{L}+C_{S}^{2}\left(  \frac{\delta}%
{L}\right)  ^{2}\frac{U^{3}}{L}%
\]

\end{theorem}

\subsection{Improving the constant multiplier}

The multiplier "$3$" of $\frac{U^{3}}{L}$\ can be reduced to $1/(1-\alpha)$
for any $\alpha>0$\ , arbitrarily close to "$1$", at the cost of a multipliers
$\alpha^{-1}$ and $\alpha^{-2}$ of the other two terms on the RHS by inserting
parameters at various points in the argument. The following is the precise result.

\begin{theorem}
$\left\langle \varepsilon_{S}\right\rangle $ satisfies: for any $0<\alpha<1,$%
\[
\left\langle \varepsilon_{S}(u)\right\rangle \leq\frac{1}{1-\alpha}\frac
{U^{3}}{L}+\frac{1}{4\alpha(1-\alpha)}\operatorname{Re}^{-1}\frac{U^{3}}%
{L}+\frac{4}{27(1-\alpha)\alpha^{2}}C_{S}^{2}\left(  \frac{\delta}{L}\right)
^{2}\frac{U^{3}}{L}.
\]

\end{theorem}

\subsection{Related work}

The energy dissipation rate is a fundamental statistic in experimental and
theoretical studies of turbulence, e.g., Sreenivasan \cite{S84}, Pope
\cite{Pope}, Frisch \cite{Frisch}, Lesieur \cite{Les97}. \ In 1968, Saffman
\cite{S68}, addressing the estimate of energy dissipation rates, $\left\langle
\varepsilon\right\rangle \simeq U^{3}/L$ , wrote\ that

\begin{center}
\textit{"This result is fundamental to an understanding of turbulence and yet
still lacks theoretical support." }- P.G. Saffman 1968
\end{center}

In 1992 Constantin and Doering \cite{CD92} made a fundamental breakthrough,
establishing a direct link between the phenomenology of energy dissipation and
that predicted for general weak solutions of shear flows directly from the
NSE. This work builds on Busse \cite{B78}, Howard \cite{H72} (and others) and
has developed in many important directions. It has been extended to shear
flows in Childress, Kerswell and Gilbert \cite{CKG01}, Kerswell \cite{K98} and
Wang \cite{Wang97}. For flows driven by body forces extensions include Doering
and Foias \cite{DF02}, Cheskidov, Doering and Petrov \cite{CDP06}\ (fractal
body forces), and \cite{L07} (helicity dissipation). The energy dissipation
rates of discretized (so the smallest scale is limited by the mesh width and
time step) flow equations was studied in \cite{JLM07}, \cite{JLK14}. Energy
dissipation in models and regularizations studied in \cite{L02}, \cite{L07},
\cite{LRS10}, \cite{LST10}. Most recently, the time averaged energy
dissipation in statistical fluctuations has led to new models and a proof of
the Boussinesq conjecture in \cite{JLK14},\cite{JL15}.

\section{The proof}

Let $||\cdot||,(\cdot,\cdot)$\ denote the usual $L^{2}(\Omega)$ norm and inner
product. Other norms are explicitly indicated by a subscript. With $|\Omega|$
the volume of the flow domain, the scale of the body force, large scale
velocity and length, $F,U,L$, are defined by%
\begin{align}
F  &  =\left(  \frac{1}{|\Omega|}||f||^{2}\right)  ^{\frac{1}{2}}\text{,
}\label{eq:ULscales}\\
U  &  =\left\langle \frac{1}{|\Omega|}||u||^{2}\right\rangle ^{\frac{1}{2}%
}\text{, }\nonumber\\
L  &  =\min\{|\Omega|^{\frac{1}{3}},\frac{F}{||\nabla f||_{L^{\infty}}}%
,\frac{F}{(\frac{1}{|\Omega|}||\nabla f||^{2})^{\frac{1}{2}}},\frac{F}%
{(\frac{1}{|\Omega|}||\nabla f||_{3}^{3})^{\frac{1}{3}}}\}.\nonumber
\end{align}
It is easy to check that $L$ has units of length and satisfies the
inequalities:%
\begin{equation}
\left.
\begin{array}
[c]{c}%
||\nabla f||_{L^{\infty}}\leq\frac{F}{L},\\
\frac{1}{|\Omega|}\int_{\Omega}|\nabla f(x)|^{2}dx\leq\frac{F^{2}}{L^{2}}\\
\frac{1}{|\Omega|}\int_{\Omega}|\nabla f(x)|^{3}dx\leq\frac{F^{3}}{L^{3}}.
\end{array}
\right\}  \label{eq:FandLproperties}%
\end{equation}
The proof is a synthesis of the model's energy balance
(\ref{eq:EnergyEquality}), the breakthrough arguments of Doering and Foias
\cite{DF02} from the NSE case with careful treatment of the Smagorinsky term.

Solutions to the Smagorinsky / Ladyzhenskaya model are known, e.g.,
\cite{C98}, \cite{DG91}, \cite{Lad69}, \cite{Lad2}, \cite{P92} \cite{G89}, to
be unique strong solutions and satisfy the energy equality%
\begin{equation}
\frac{1}{2|\Omega|}||u(T)||^{2}+\int_{0}^{T}\varepsilon_{S}(u)dt=\frac
{1}{2|\Omega|}||u_{0}||^{2}+\int_{0}^{T}\frac{1}{|\Omega|}(f,u(t))\,dt.
\label{eq:EnergyEquality}%
\end{equation}
Here $\varepsilon_{S}(u)$ = $\varepsilon_{0}(u)+\varepsilon_{\delta}(u),$
where%
\begin{align*}
\varepsilon_{0}(u)  &  :=\frac{\nu}{|\Omega|}||\nabla u(t)||^{2}\text{ and }\\
\varepsilon_{\delta}(u)  &  :=\frac{(C_{S}\delta)^{2}}{|\Omega|}||\nabla
u(t)||_{L^{3}}^{3}.
\end{align*}
From (\ref{eq:EnergyEquality}) and standard arguments it follows that%
\begin{align}
\sup_{t\in(0,\infty)}||u(t)||^{2}  &  \leq C(data)<\infty\text{ and}%
\label{eq:Boundsonuand gradu}\\
\text{ }\frac{1}{T}\int_{0}^{T}\varepsilon_{S}(u)dt  &  \leq C(data)<\infty
.\nonumber
\end{align}
Averaging (\ref{eq:EnergyEquality}) over $[0,T]$, applying the
Cauchy-Schwarz\ inequality in time and (\ref{eq:Boundsonuand gradu}) yields
\begin{align}
\frac{1}{T}\int_{0}^{T}\varepsilon_{S}(u)dt  &  =\mathcal{O}(\frac{1}%
{T})+\frac{1}{T}\int_{0}^{T}\frac{1}{|\Omega|}(f,u(t))\,dt\nonumber\\
&  \leq\mathcal{O}(\frac{1}{T})+F\left(  \frac{1}{T}\int_{0}^{T}\frac
{1}{|\Omega|}||u||^{2}dt\right)  ^{\frac{1}{2}}. \label{eq:Step1}%
\end{align}
To bound the RHS, take the inner product of (\ref{eq:SmagorinskyModel}) with
$f(x)$ , integrate by parts and average over $[0,T]$. This gives%
\begin{gather}
F^{2}=\frac{(u(T)-u_{0},f)}{T|\Omega|}-\frac{1}{T}\int_{0}^{T}\frac{1}%
{|\Omega|}(uu,\nabla f)dt+\label{eq:Step2}\\
+\frac{1}{T}\int_{0}^{T}\frac{\nu}{|\Omega|}(\nabla u,\nabla f)dt+\frac{1}%
{T}\int_{0}^{T}\frac{1}{|\Omega|}\left(  \left(  C_{S}\delta\right)
^{2}|\nabla u|\nabla u,\nabla f\right)  dt\,.\nonumber
\end{gather}
Of the four terms on the last RHS, by (\ref{eq:Boundsonuand gradu}) the first
term is\ $\mathcal{O}(1/T)$. The second and third terms are bounded using the
Cauchy-Schwarz-Young inequality and (\ref{eq:FandLproperties}) by%
\begin{align*}
\left\vert \frac{1}{T|\Omega|}\int_{0}^{T}(uu,\nabla f)dt\right\vert  &
\leq||\nabla f||_{L^{\infty}}\frac{1}{T|\Omega|}\int_{0}^{T}||u||^{2}dt\\
&  \leq\frac{F}{L}\left(  \frac{1}{T}\int_{0}^{T}\frac{1}{|\Omega|}%
||u||^{2}dt\right)  ,\text{ }\\
\left\vert \frac{1}{T}\int_{0}^{T}\frac{\nu}{|\Omega|}(\nabla u,\nabla
f)dt\right\vert  &  \leq\left(  \frac{1}{T}\int_{0}^{T}\frac{\nu^{2}}%
{|\Omega|}||\nabla u||^{2}dt\right)  ^{\frac{1}{2}}\left(  \frac{1}{T}\int%
_{0}^{T}\frac{1}{|\Omega|}||\nabla f||^{2}dt\right)  ^{\frac{1}{2}}\\
&  \leq\left(  \frac{1}{T}\int_{0}^{T}\varepsilon_{0}(u)dt\right)  ^{\frac
{1}{2}}\sqrt{\nu}\frac{F}{L}\\
&  \leq\frac{2}{3}U^{-1}F\frac{1}{T}\int_{0}^{T}\varepsilon_{0}(u)dt+\frac
{3}{8}UF\frac{\nu}{L^{2}}.
\end{align*}
The fourth, Smagorinsky, term, is estimated using H\"{o}lder's inequality as
follows%
\begin{align*}
\left\vert \frac{1}{T|\Omega|}\int_{0}^{T}\left(  \left(  C_{S}\delta\right)
^{2}|\nabla u|\nabla u,\nabla f\right)  dt\right\vert  &  \leq\frac{\left(
C_{S}\delta\right)  ^{2}}{|\Omega|}\frac{1}{T}\int_{0}^{T}||\nabla u||_{L^{3}%
}^{2}dt||\nabla f||_{L^{3}}\\
&  \leq\frac{F}{L}\left(  C_{S}\delta\right)  ^{2/3}\frac{1}{T}\int_{0}%
^{T}\varepsilon_{\delta}(u)^{\frac{2}{3}}dt.
\end{align*}
Insert multipliers of $U^{2/3}$ and~$U^{-2/3}$ \ in the two terms.
Using\footnote{More generally, for conjugate ($1/p+1/q=1$) exponents and
$a>0,b>0$: $ab\leq(1/p)a^{p}+(1/q)b^{q\text{ }}.$}
\[
ab\leq(2/3)a^{3/2}+(1/3)b^{3\text{ }}%
\]
(conjugate exponents $3/2$ and $3$) gives%
\[
\frac{1}{T}\int_{0}^{T}\left(  \frac{U^{2/3}}{L}\left(  C_{S}\delta\right)
^{2/3}\right)  \left(  U^{-2/3}\varepsilon_{\delta}(u)^{\frac{2}{3}}\right)
dt\leq\frac{2}{3}\frac{1}{U}\frac{1}{T}\int_{0}^{T}\varepsilon_{\delta
}(u)dt+\frac{U^{2}}{3}\frac{\left(  C_{S}\delta\right)  ^{2}}{L^{3}}.
\]
Using these three estimates in (\ref{eq:Step2}) yields%
\begin{align*}
F  &  \leq\mathcal{O}(\frac{1}{T})+\frac{1}{L}\frac{1}{T}\int_{0}^{T}\frac
{1}{|\Omega|}||u||^{2}dt+\\
&  +\frac{3}{8}\frac{U\nu}{L^{2}}+\frac{2}{3}\frac{1}{U}\frac{1}{T}\int%
_{0}^{T}\varepsilon_{S}(u)dt+\frac{U^{2}}{3}\frac{\left(  C_{S}\delta\right)
^{2}}{L^{3}}\,.
\end{align*}
Using this estimate for $F$ in (\ref{eq:Step1}) \ gives%
\begin{gather*}
\frac{1}{T}\int_{0}^{T}\varepsilon_{S}(u)dt\leq\mathcal{O}(\frac{1}%
{T})+F\left(  \frac{1}{T}\int_{0}^{T}\frac{1}{|\Omega|}||u||^{2}dt\right)
^{\frac{1}{2}}\\
\leq\mathcal{O}(\frac{1}{T})+\left(  \frac{1}{T|\Omega|}\int_{0}^{T}%
||u||^{2}dt\right)  ^{\frac{1}{2}}\times\\
\left(  \frac{1}{LT|\Omega|}\int_{0}^{T}||u||^{2}dt+\frac{3}{8}\frac{U\nu
}{L^{2}}+\frac{2}{3}\frac{1}{UT}\int_{0}^{T}\varepsilon_{S}(u)dt+\frac
{U^{2}\left(  C_{S}\delta\right)  ^{2}}{3L^{3}}\right)  \,\,.
\end{gather*}
Taking the limit superior, which exists by\ (\ref{eq:Boundsonuand gradu}), as
$T\rightarrow\infty$ we obtain%
\[
\left\langle \varepsilon_{S}(u)\right\rangle \leq\frac{U^{3}}{L}+\frac{3}%
{8}\frac{U^{3}}{L}\frac{\nu}{LU}+\frac{2}{3}\left\langle \varepsilon
_{S}(u)\right\rangle +\frac{U^{3}}{3L}\frac{\left(  C_{S}\delta\right)  ^{2}%
}{L^{2}}.
\]
Thus, as claimed,%
\[
\left\langle \varepsilon_{S}(u)\right\rangle \leq3\frac{U^{3}}{L}+\frac{9}%
{8}Re^{-1}\frac{U^{3}}{L}+\frac{U^{3}}{L}C_{S}^{2}\left(  \frac{\delta}%
{L}\right)  ^{2}.
\]

\section{Conclusions for the Smagorinsky Model}

Comparing the estimate
\[
\left\langle \varepsilon_{S}\right\rangle \simeq\frac{U^{3}}{L}%
\]
herein for periodic boundary conditions with
\[
\left\langle \varepsilon_{S}\right\rangle \simeq\lbrack1+C_{S}^{2}\left(
\frac{\delta}{L}\right)  ^{2}(1+Re^{2})]\frac{U^{3}}{L}%
\]
\ in \cite{L07} for shear flows with boundary layers strongly suggests the
often reported model over dissipation is

\begin{center}
\textit{due to the action of the model viscosity in boundary layers}
\end{center}

rather than in interior small scales generated by the turbulent cascade.
Practice addresses this over dissipation with damping functions (e.g., van
Driest damping \cite{vD},\cite{Sagaut}), modelled boundary conditions called
near wall models, e.g., \cite{PB02}, \cite{JLS04}, \cite{JL06}, restriction of
the model induced dissipation to the smallest resolved scales \cite{HOM01} and
Germano's dynamic (self-adaptive) selection of $C_{S}=C_{S}(x,t)$ ,
\cite{GPMC91}, that also reduces $C_{S}$\ near walls. Thus, analysis of
Smagorinsky\ model\ dissipation for\ shear flows \textit{including these
modifications} is therefore an important open problem.

\end{document}